# A Degree Condition for Dominating Cycles in $t$-tough Graphs with $t > 1$

Zh.G. Nikoghosyan[*]

January 7, 2012

**Abstract**

Let $G$ be a $t$-tough graph of order $n$ and minimum degree $\delta$ with $t > 1$. It is proved that if $\delta \geq (n-2)/3$ then each longest cycle in $G$ is a dominating cycle.

## 1 Introduction

Only finite undirected graphs without loops or multiple edges are considered. We reserve $n$, $\delta$, $\kappa$ and $\tau$ to denote the number of vertices (order), the minimum degree, connectivity and the toughness of a graph, respectively. A good reference for any undefined terms is [3].

The earliest sufficient condition for dominating cycles was developed in 1971 due to Nash-Williams [4].

**Theorem A [4].** Let $G$ be a 2-connected graph. If $\delta \geq \frac{1}{3}(n+2)$ then each longest cycle in $G$ is a dominating cycle.

In 1979, Bigalke and Jung [2] proved that the minimum degree condition $\delta \geq (n+2)/3$ in Theorem A can be slightly relaxed by replacing the 2-connectivity condition with stronger 1-tough condition.

**Theorem B [2].** Let $G$ be a 1-tough graph. If $\delta \geq \frac{1}{3}n$ then each longest cycle in $G$ is a dominating cycle.

In this paper we prove that the minimum degree bound $n/3$ in Theorem B can be lowered to $\delta \geq (n-2)/3$ when $\tau > 1$.

---

[*]G.G. Nicoghossian (up to 1997)

**Theorem 1**. Let $G$ be a graph with $\tau > 1$. If $\delta \geq \frac{1}{3}(n-2)$ then each longest cycle in $G$ is a dominating cycle.

To prove Theorem 1, we need two known lower bounds for the circumference, the length of a longest cycle in a graph. The first bound is due to Voss and Zuluaga [5] concerning the alternative existence of long cycles and dominating cycles in 3-connected graphs.

**Theorem C [5]**. Let $G$ be a 3-connected graph. Then either $G$ has a cycle of length at least $3\delta - 3$ or each longest cycle in $G$ is a dominating cycle.

The second bound is due to Bauer and Schmeichel [1] concerning the alternative existence of long cycles and Hamilton cycles in 1-tough graphs.

**Theorem D [1]**. Let $G$ be a graph with $\tau \geq 1$. Then either $G$ has a Hamilton cycle or has a cycle of length at least $2\delta + 2$.

## 2 Notations and preliminaries

The set of vertices of a graph $G$ is denoted by $V(G)$ and the set of edges by $E(G)$. For $S$ a subset of $V(G)$, we denote by $G\backslash S$ the maximum subgraph of $G$ with vertex set $V(G)\backslash S$. We write $G[S]$ for the subgraph of $G$ induced by $S$. For a subgraph $H$ of $G$ we use $G\backslash H$ short for $G\backslash V(H)$. The neighborhood of a vertex $x \in V(G)$ will be denoted by $N(x)$. Furthermore, for a subgraph $H$ of $G$ and $x \in V(G)$, we define $N_H(x) = N(x) \cap V(H)$ and $d_H(x) = |N_H(x)|$. Let $s(G)$ denote the number of components of a graph $G$. A graph $G$ is $t$-tough if $|S| \geq ts(G\backslash S)$ for every subset $S$ of the vertex set $V(G)$ with $s(G\backslash S) > 1$. The toughness of $G$, denoted $\tau(G)$, is the maximum value of $t$ for which $G$ is $t$-tough (taking $\tau(K_n) = \infty$ for all $n \geq 1$).

A simple cycle (or just a cycle) $C$ of length $t$ is a sequence $v_1v_2...v_tv_1$ of distinct vertices $v_1,...,v_t$ with $v_iv_{i+1} \in E(G)$ for each $i \in \{1,...,t\}$, where $v_{t+1} = v_1$. When $t = 2$, the cycle $C = v_1v_2v_1$ on two vertices $v_1, v_2$ coincides with the edge $v_1v_2$, and when $t = 1$, the cycle $C = v_1$ coincides with the vertex $v_1$. So, all vertices and edges in a graph can be considered as cycles of lengths 1 and 2, respectively. A graph $G$ is hamiltonian if $G$ contains a Hamilton cycle, i.e. a cycle of length $n$. A cycle $C$ in $G$ is dominating if $G\backslash C$ is edgeless.

Paths and cycles in a graph $G$ are considered as subgraphs of $G$. If $Q$ is a path or a cycle, then the length of $Q$, denoted by $|Q|$, is $|E(Q)|$. We write $Q$ with a given orientation by $\overrightarrow{Q}$. For $x, y \in V(Q)$, we denote by $x\overrightarrow{Q}y$ the subpath of $Q$ in the chosen direction from $x$ to $y$.

**Special definitions**. Let $G$ be a graph, $C$ a longest cycle in $G$ and $P = x\overrightarrow{P}y$ a longest path in $G\backslash C$ of length $\overline{p} \geq 0$. Let $\xi_1, \xi_2, ..., \xi_s$ be the elements of

$N_C(x) \cup N_C(y)$ occuring on $C$ in a consecutive order. Set

$$I_i = \xi_i \overrightarrow{C} \xi_{i+1}, \ I_i^* = \xi_i^+ \overrightarrow{C} \xi_{i+1}^- \ \ (i = 1, 2, ..., s),$$

where $\xi_{s+1} = \xi_1$.

(∗1) We call $I_1, I_2, ..., I_s$ elementary segments on $C$ created by $N_C(x) \cup N_C(y)$.

(∗2) We call a path $L = z \overrightarrow{L} w$ an intermediate path between two distinct elementary segments $I_a$ and $I_b$ if

$$z \in V(I_a^*), \ w \in V(I_b^*), \ V(L) \cap V(C \cup P) = \{z, w\}.$$

(∗3) The set of all intermediate paths between elementary segments $I_{i_1}, I_{i_2}, ..., I_{i_t}$ will be denoted by $\Upsilon(I_{i_1}, I_{i_2}, ..., I_{i_t})$ .

**Lemma 1.** Let $G$ be a graph, $C$ a longest cycle in $G$ and $P = x \overrightarrow{P} y$ a longest path in $G \backslash C$ of length $\overline{p} \geq 1$. If $|N_C(x)| \geq 2$, $|N_C(y)| \geq 2$ and $N_C(x) \neq N_C(y)$ then

$$|C| \geq \begin{cases} 3\delta + \max\{\sigma_1, \sigma_2\} - 1 \geq 3\delta & \text{if} \quad \overline{p} = 1, \\ \max\{2\overline{p} + 8, 4\delta - 2\overline{p}\} & \text{if} \quad \overline{p} \geq 2, \end{cases}$$

where $\sigma_1 = |N_C(x) \backslash N_C(y)|$ and $\sigma_2 = |N_C(y) \backslash N_C(x)|$.

**Lemma 2.** Let $G$ be a graph, $C$ a longest cycle in $G$ and $P = x \overrightarrow{P} y$ a longest path in $G \backslash C$ of length $\overline{p} \geq 0$. If $N_C(x) = N_C(y)$ and $|N_C(x)| \geq 2$ then for each elementary segments $I_a$ and $I_b$ induced by $N_C(x) \cup N_C(y)$,

(a1) if $L$ is an intermediate path between $I_a$ and $I_b$ then

$$|I_a| + |I_b| \geq 2\overline{p} + 2|L| + 4,$$

(a2) if $\Upsilon(I_a, I_b) \subseteq E(G)$ and $|\Upsilon(I_a, I_b)| = i$ for some $i \in \{1, 2, 3\}$ then

$$|I_a| + |I_b| \geq 2\overline{p} + i + 5,$$

(a3) if $\Upsilon(I_a, I_b) \subseteq E(G)$ and $\Upsilon(I_a, I_b)$ contains two independent intermediate edges then

$$|I_a| + |I_b| \geq 2\overline{p} + 8.$$

## 3 Proofs

**Proof of Lemma 1**. Put

$$A_1 = N_C(x) \backslash N_C(y), \ A_2 = N_C(y) \backslash N_C(x), \ M = N_C(x) \cap N_C(y).$$

By the hypothesis, $N_C(x) \neq N_C(y)$, implying that

$$\max\{|A_1|, |A_2|\} \geq 1.$$

Let $\xi_1, \xi_2, ..., \xi_s$ be the elements of $N_C(x) \cup N_C(y)$ occuring on $C$ in a consecutive order. Put $I_i = \xi_i \overrightarrow{C} \xi_{i+1}$ $(i = 1, 2, ..., s)$, where $\xi_{s+1} = \xi_1$. Clearly, $s = |A_1| + |A_2| + |M|$. Since $C$ is extreme, $|I_i| \geq 2$ $(i = 1, 2, ..., s)$. Next, if $\{\xi_i, \xi_{i+1}\} \cap M \neq \emptyset$ for some $i \in \{1, 2, ..., s\}$ then $|I_i| \geq \overline{p} + 2$. Further, if either $\xi_i \in A_1$, $\xi_{i+1} \in A_2$ or $\xi_i \in A_2$, $\xi_{i+1} \in A_1$ then again $|I_i| \geq \overline{p} + 2$.

**Case 1**. $\overline{p} = 1$.
**Case 1.1**. $|A_i| \geq 1$ $(i = 1, 2)$.

It follows that among $I_1, I_2, ..., I_s$ there are $|M| + 2$ segments of length at least $\overline{p} + 2$. Observing also that each of the remaining $s - (|M| + 2)$ segments has a length at least 2, we have

$$|C| \geq (\overline{p} + 2)(|M| + 2) + 2(s - |M| - 2)$$

$$= 3(|M| + 2) + 2(|A_1| + |A_2| - 2)$$

$$= 2|A_1| + 2|A_2| + 3|M| + 2.$$

Since $|A_1| = d(x) - |M| - 1$ and $|A_2| = d(y) - |M| - 1$,

$$|C| \geq 2d(x) + 2d(y) - |M| - 2 \geq 3\delta + d(x) - |M| - 2.$$

Recalling that $d(x) = |M| + |A_1| + 1$, we get

$$|C| \geq 3\delta + |A_1| - 1 = 3\delta + \sigma_1 - 1.$$

Analogously, $|C| \geq 3\delta + \sigma_2 - 1$. So,

$$|C| \geq 3\delta + \max\{\sigma_1, \sigma_2\} - 1 \geq 3\delta.$$

**Case 1.2**. Either $|A_1| \geq 1, |A_2| = 0$ or $|A_1| = 0, |A_2| \geq 1$.

Assume w.l.o.g. that $|A_1| \geq 1$ and $|A_2| = 0$, i.e. $|N_C(y)| = |M| \geq 2$ and $s = |A_1| + |M|$. Hence, among $I_1, I_2, ..., I_s$ there are $|M| + 1$ segments of length at least $\overline{p} + 2 = 3$. Taking into account that each of the remaining $s - (|M| + 1)$ segments has a length at least 2 and $|M| + 1 = d(y)$, we get

$$|C| \geq 3(|M| + 1) + 2(s - |M| - 1) = 3d(y) + 2(|A_1| - 1)$$

$$\geq 3\delta + |A_1| - 1 = 3\delta + \max\{\sigma_1, \sigma_2\} - 1 \geq 3\delta.$$

**Case 2**. $\overline{p} \geq 2$.

We first prove that $|C| \geq 2\overline{p} + 8$. Since $|N_C(x)| \geq 2$ and $|N_C(y)| \geq 2$, there are at least two segments among $I_1, I_2, ..., I_s$ of length at least $\overline{p} + 2$. If $|M| = 0$ then clearly $s \geq 4$ and

$$|C| \geq 2(\overline{p} + 2) + 2(s - 2) \geq 2\overline{p} + 8.$$

Otherwise, since $\max\{|A_1|, |A_2|\} \geq 1$, there are at least three elementary segments of length at least $\overline{p} + 2$, that is

$$|C| \geq 3(\overline{p} + 2) \geq 2\overline{p} + 8.$$

So, in any case, $|C| \geq 2\overline{p} + 8$.

To prove that $|C| \geq 4\delta - 2\overline{p}$, we distinguish two main cases.

**Case 2.1**. $|A_i| \geq 1$ $(i = 1, 2)$.

It follows that among $I_1, I_2, ..., I_s$ there are $|M| + 2$ segments of length at least $\overline{p} + 2$. Further, since each of the remaining $s - (|M| + 2)$ segments has a length at least 2, we get

$$|C| \geq (\overline{p} + 2)(|M| + 2) + 2(s - |M| - 2)$$
$$= (\overline{p} - 2)|M| + (2\overline{p} + 4|M| + 4) + 2(|A_1| + |A_2| - 2)$$
$$\geq 2|A_1| + 2|A_2| + 4|M| + 2\overline{p}.$$

Observing also that

$$|A_1| + |M| + \overline{p} \geq d(x), \quad |A_2| + |M| + \overline{p} \geq d(y),$$

we have

$$2|A_1| + 2|A_2| + 4|M| + 2\overline{p}$$
$$\geq 2d(x) + 2d(y) - 2\overline{p} \geq 4\delta - 2\overline{p},$$

implying that $|C| \geq 4\delta - 2\overline{p}$.

**Case 2.2**. Either $|A_1| \geq 1, |A_2| = 0$ or $|A_1| = 0, |A_2| \geq 1$.

Assume w.l.o.g. that $|A_1| \geq 1$ and $|A_2| = 0$, i.e. $|N_C(y)| = |M| \geq 2$ and $s = |A_1| + |M|$. It follows that among $I_1, I_2, ..., I_s$ there are $|M| + 1$ segments of length at least $\overline{p} + 2$. Observing also that $|M| + \overline{p} \geq d(y) \geq \delta$, i.e. $2\overline{p} + 4|M| \geq 4\delta - 2\overline{p}$, we get

$$|C| \geq (\overline{p} + 2)(|M| + 1) \geq (\overline{p} - 2)(|M| - 1) + 2\overline{p} + 4|M|$$
$$\geq 2\overline{p} + 4|M| \geq 4\delta - 2\overline{p}. \qquad \blacksquare$$

**Proof of Lemma 2**. Let $\xi_1, \xi_2, ..., \xi_s$ be the elements of $N_C(x)$ occuring on $C$ in a consecutive order. Put $I_i = \xi_i \overrightarrow{C} \xi_{i+1}$ $(i = 1, 2, ..., s)$, where $\xi_{s+1} = \xi_1$. To prove (a1), let $L = z\overrightarrow{L}w$ be an intermediate path between elementary segments $I_a$ and $I_b$ with $z \in V(I_a^*)$ and $w \in V(I_b^*)$. Put

$$|\xi_a \overrightarrow{C} z| = d_1, \ |z\overrightarrow{C}\xi_{a+1}| = d_2, \ |\xi_b \overrightarrow{C} w| = d_3, \ |w\overrightarrow{C}\xi_{b+1}| = d_4,$$
$$C' = \xi_a x \overrightarrow{P} y \xi_b \overleftarrow{C} z \overrightarrow{L} w \overrightarrow{C} \xi_a.$$

Clearly,

$$|C'| = |C| - d_1 - d_3 + |L| + |P| + 2.$$

Since $C$ is extreme, we have $|C| \geq |C'|$, implying that $d_1 + d_3 \geq \overline{p} + |L| + 2$. By a symmetric argument, $d_2 + d_4 \geq \overline{p} + |L| + 2$. Hence

$$|I_a| + |I_b| = \sum_{i=1}^{4} d_i \geq 2\overline{p} + 2|L| + 4.$$

The proof of $(a1)$ is complete. To proof $(a2)$ and $(a3)$, let $\Upsilon(I_a, I_b) \subseteq E(G)$ and $|\Upsilon(I_a, I_b)| = i$ for some $i \in \{1, 2, 3\}$.

**Case 1.** $i = 1$.

It follows that $\Upsilon(I_a, I_b)$ consists of a unique intermediate edge $L = zw$. By $(a1)$,
$$|I_a| + |I_b| \geq 2\overline{p} + 2|L| + 4 = 2\overline{p} + 6.$$

**Case 2.** $i = 2$.

It follows that $\Upsilon(I_a, I_b)$ consists of two edges $e_1, e_2$. Put $e_1 = z_1 w_1$ and $e_2 = z_2 w_2$, where $\{z_1, z_2\} \subseteq V(I_a^*)$ and $\{w_1, w_2\} \subseteq V(I_b^*)$.

**Case 2.1.** $z_1 \neq z_2$ and $w_1 \neq w_2$.

Assume w.l.o.g. that $z_1$ and $z_2$ occur in this order on $I_a$.

**Case 2.1.1.** $w_2$ and $w_1$ occur in this order on $I_b$.

Put
$$|\xi_a \overrightarrow{C} z_1| = d_1, \ |z_1 \overrightarrow{C} z_2| = d_2, \ |z_2 \overrightarrow{C} \xi_{a+1}| = d_3,$$
$$|\xi_b \overrightarrow{C} w_2| = d_4, \ |w_2 \overrightarrow{C} w_1| = d_5, \ |w_1 \overrightarrow{C} \xi_{b+1}| = d_6,$$
$$C' = \xi_a \overrightarrow{C} z_1 w_1 \overleftarrow{C} w_2 z_2 \overrightarrow{C} \xi_b x \overrightarrow{P} y \xi_{b+1} \overrightarrow{C} \xi_a.$$

Clearly,
$$|C'| = |C| - d_2 - d_4 - d_6 + |\{e_1\}| + |\{e_2\}| + |P| + 2$$
$$= |C| - d_2 - d_4 - d_6 + \overline{p} + 4.$$

Since $C$ is extreme, $|C| \geq |C'|$, implying that $d_2 + d_4 + d_6 \geq \overline{p} + 4$. By a symmetric argument, $d_1 + d_3 + d_5 \geq \overline{p} + 4$. Hence
$$|I_a| + |I_b| = \sum_{i=1}^{6} d_i \geq 2\overline{p} + 8.$$

**Case 2.1.2.** $w_1$ and $w_2$ occur in this order on $I_b$.

Putting
$$C' = \xi_a \overrightarrow{C} z_1 w_1 \overrightarrow{C} w_2 z_2 \overrightarrow{C} \xi_b x \overrightarrow{P} y \xi_{b+1} \overrightarrow{C} \xi_a,$$
we can argue as in Case 2.1.1.

**Case 2.2.** Either $z_1 = z_2$, $w_1 \neq w_2$ or $z_1 \neq z_2$, $w_1 = w_2$.

Assume w.l.o.g. that $z_1 \neq z_2$, $w_1 = w_2$ and $z_1, z_2$ occur in this order on $I_a$.
Put
$$|\xi_a \overrightarrow{C} z_1| = d_1, \ |z_1 \overrightarrow{C} z_2| = d_2, \ |z_2 \overrightarrow{C} \xi_{a+1}| = d_3,$$
$$|\xi_b \overrightarrow{C} w_1| = d_4, \ |w_1 \overrightarrow{C} \xi_{b+1}| = d_5,$$
$$C' = \xi_a x \overrightarrow{P} y \xi_b \overleftarrow{C} z_1 w_1 \overrightarrow{C} \xi_a,$$
$$C'' = \xi_a \overrightarrow{C} z_2 w_1 \overleftarrow{C} \xi_{a+1} x \overrightarrow{P} y \xi_{b+1} \overrightarrow{C} \xi_a.$$

Clearly,
$$|C'| = |C| - d_1 - d_4 + |\{e_1\}| + |P| + 2 = |C| - d_1 - d_4 + \overline{p} + 3,$$
$$|C''| = |C| - d_3 - d_5 + |\{e_2\}| + |P| + 2 = |C| - d_3 - d_5 + \overline{p} + 3.$$

Since $C$ is extreme, $|C| \geq |C'|$ and $|C| \geq |C''|$, implying that
$$d_1 + d_4 \geq \overline{p} + 3, \ d_3 + d_5 \geq \overline{p} + 3.$$

Hence,
$$|I_a| + |I_b| = \sum_{i=1}^{5} d_i \geq d_1 + d_3 + d_4 + d_5 + 1 \geq 2\overline{p} + 7.$$

**Case 3.** $i = 3$.

It follows that $\Upsilon(I_a, I_b)$ consists of three edges $e_1, e_2, e_3$. Let $e_i = z_i w_i$ ($i = 1, 2, 3$), where $\{z_1, z_2, z_3\} \subseteq V(I_a^*)$ and $\{w_1, w_2, w_3\} \subseteq V(I_b^*)$. If there are two independent edges among $e_1, e_2, e_3$ then we can argue as in Case 2.1. Otherwise, we can assume w.l.o.g. that $w_1 = w_2 = w_3$ and $z_1, z_2, z_3$ occur in this order on $I_a$. Put
$$|\xi_a \overrightarrow{C} z_1| = d_1, \ |z_1 \overrightarrow{C} z_2| = d_2, \ |z_2 \overrightarrow{C} z_3| = d_3,$$
$$|z_3 \overrightarrow{C} \xi_{a+1}| = d_4, \ |\xi_b \overrightarrow{C} w_1| = d_5, \ |w_1 \overrightarrow{C} \xi_{b+1}| = d_6,$$
$$C' = \xi_a x \overrightarrow{P} y \xi_b \overleftarrow{C} z_1 w_1 \overrightarrow{C} \xi_a,$$
$$C'' = \xi_a \overrightarrow{C} z_3 w_1 \overleftarrow{C} \xi_{a+1} x \overrightarrow{P} y \xi_{b+1} \overrightarrow{C} \xi_a.$$

Clearly,
$$|C'| = |C| - d_1 - d_5 + |\{e_1\}| + \overline{p} + 2,$$
$$|C''| = |C| - d_4 - d_6 + |\{e_3\}| + \overline{p} + 2.$$

Since $C$ is extreme, we have $|C| \geq |C'|$ and $|C| \geq |C''|$, implying that
$$d_1 + d_5 \geq \overline{p} + 3, \ d_4 + d_6 \geq \overline{p} + 3.$$

Hence,
$$|I_a| + |I_b| = \sum_{i=1}^{6} d_i \geq d_1 + d_4 + d_5 + d_6 + 2 \geq 2\overline{p} + 8. \ \blacksquare$$

**Proof of Theorem 1.** Let $C$ be a longest cycle in $G$ and $P = x_1 \overrightarrow{P} x_2$ a longest path in $G \backslash C$ of length $\overline{p}$. If $|V(P)| \leq 1$ then $C$ is a dominating cycle and we are done. Let $|V(P)| \geq 2$, that is $\overline{p} \geq 1$. If $\kappa \leq 2$ then clearly $\tau \leq 1$, contradicting the hypothesis. Let $\kappa \geq 3$. By Theorem C, $|C| \geq 3\delta - 3$. On the other hand, by the hypothesis, $|C| + \overline{p} + 1 \leq n \leq 3\delta + 2$, implying that
$$3\delta - 3 \leq |C| \leq 3\delta - \overline{p} + 1. \tag{1}$$

Let $\xi_1, \xi_2, ..., \xi_s$ be the elements of $N_C(x_1) \cup N_C(x_2)$ occuring on $C$ in a consecutive order. Put

$$I_i = \xi_i \overrightarrow{C} \xi_{i+1}, \ I_i^* = \xi_i^+ \overrightarrow{C} \xi_{i+1}^- \ (i = 1, 2, ..., s),$$

where $\xi_{s+1} = \xi_1$.

**Claim 1.** Let $N_C(y) \subseteq \{\xi_1, \xi_2, ..., \xi_s\}$ for each $y \in V(P)$. If $\Upsilon(I_1, I_2, ..., I_s) \subseteq E(G)$ and the edges in $\Upsilon(I_1, I_2, ..., I_s)$ form a star then $\tau \leq 1$.

**Proof of Claim 1.** By the hypothesis, there is a vertex $z$ belonging to all edges in $\Upsilon(I_1, I_2, ..., I_s)$. Then $G \backslash \{\xi_1, \xi_2, ..., \xi_s, z\}$ has at least $s + 1$ connected components, that is $\tau \leq 1$. $\Delta$

**Case 1.** $\overline{p} = 1$.
By (1),
$$3\delta - 3 \leq |C| \leq 3\delta. \tag{2}$$

Since $\delta \geq \kappa \geq 3$, we have $|N_C(x_i)| \geq \delta - \overline{p} = \delta - 1 \geq 2 \ (i = 1, 2)$.

**Case 1.1.** $N_C(x_1) \neq N_C(x_2)$.
It follows that $\max\{\sigma_1, \sigma_2\} \geq 1$, where

$$\sigma_1 = |N_C(x_1) \backslash N_C(x_2)|, \quad \sigma_2 = |N_C(x_2) \backslash N_C(x_1)|.$$

By Lemma 1, $|C| \geq 3\delta$. Recalling (2), we get $|C| = 3\delta$. If $\max\{\sigma_1, \sigma_2\} \geq 2$ then by Lemma 1, $|C| \geq 3\delta + 1$, contradicting (2). Let $\max\{\sigma_1, \sigma_2\} = 1$. This implies $s \geq \delta$ and $|I_i| \geq 3 \ (i = 1, 2, ..., s)$. If $s \geq \delta + 1$ then $|C| \geq 3s \geq 3\delta + 3$, again contradicting (2). Let $s = \delta$, that is $|I_i| = 3 \ (i = 1, 2, ..., s)$. By Lemma 2, $\Upsilon(I_1, I_2, ..., I_s) = \emptyset$, which by Claim 1 yields $\tau \leq 1$, contradicting the hypothesis.

**Case 1.2.** $N_C(x_1) = N_C(x_2)$.
Clearly, $s = |N_C(x_1)| \geq \delta - \overline{p} \geq \delta - 1$. If $s \geq \delta$ then $|C| \geq 3s \geq 3\delta$ and we can argue as in Case 1.1. Let $s = \delta - 1$.
From (2) and Lemma 2, we can easily obtain the following.

**Claim 2.**
(1) $|I_i| + |I_j| \leq 9$ for each distinct $i, j \in \{1, 2, ..., s\}$.
(2) If $|I_a| + |I_b| = 9$ for some distinct $a, b \in \{1, 2, ..., s\}$ then $|I_i| = 3$ for each $i \in \{1, 2, ..., s\} \backslash \{a, b\}$.
(3) If $|I_a| = 6$ for some $a \in \{1, 2, ..., s\}$ then $|I_i| = 3$ for each $i \in \{1, 2, ..., s\} \backslash \{a\}$.
(4) There are at most three segments of length at least 4.
(5) If $|I_a| \geq 4$, $|I_b| \geq 4$, $|I_c| \geq 4$ for some distinct $a, b, c \in \{1, 2, ..., s\}$ then $|I_a| = |I_b| = |I_c| = 4$.

If $\Upsilon(I_1, I_2, ..., I_s) = \emptyset$ then by Claim 1, $\tau \leq 1$, contradicting the hypothesis. Otherwise, $\Upsilon(I_a, I_b) \neq \emptyset$ for some distinct $a, b \in \{1, 2, ..., s\}$. By definition, there is an intermediate path $L$ between $I_a$ and $I_b$. If $|L| \geq 2$ then by Lemma 2,
$$|I_a| + |I_b| \geq 2\overline{p} + 2|L| + 4 \geq 10,$$
contradicting Claim 2. Otherwise, $|L| = 1$ and therefore,
$$\Upsilon(I_1, I_2, ..., I_s) \subseteq E(G).$$
By Lemma 2, $|I_a| + |I_b| \geq 2\overline{p} + 6 = 8$. Combining this with Claim 2, we have
$$8 \leq |I_a| + |I_b| \leq 9.$$
Let $L = yz$, where $y \in V(I_a^*)$ and $z \in V(I_b^*)$. Put
$$C_1 = \xi_a x_1 x_2 \xi_b \overleftarrow{C} yz \overrightarrow{C} \xi_a,$$
$$C_2 = \xi_a \overrightarrow{C} yz \overleftarrow{C} \xi_{a+1} x_1 x_2 \xi_{b+1} \overrightarrow{C} \xi_a.$$

**Case 1.2.1.** $|I_a| + |I_b| = 8$.

Since $|I_i| \geq 3$ $(i = 1, 2, ..., s)$, we can assume w.l.o.g. that either $|I_a| = 3$, $|I_b| = 5$ or $|I_a| = |I_b| = 4$.

**Case 1.2.1.1.** $|I_a| = 3$ and $|I_b| = 5$.

Put $I_a = \xi_a w_1 w_2 \xi_{a+1}$ and $I_b = \xi_b w_3 w_4 w_5 w_6 \xi_{b+1}$. Assume w.l.o.g. that $y = w_2$. If $z = w_3$ then $|C_1| > |C|$, a contradiction. Further, if $z \in \{w_5, w_6\}$ then $|C_2| > |C|$, a contradiction. Hence, $z = w_4$. To determine the possible neighborhood of $w_1$, we first observe that if $w_1 w_4 \in E(G)$ then
$$\xi_a x_1 x_2 \xi_b \overleftarrow{C} w_1 w_4 \overrightarrow{C} \xi_a$$
is longer than $C$, a contradiction. Let $w_1 w_4 \notin E(G)$. Therefore, if $N(w_1) \cap V(I_b^*) \neq \emptyset$ then there exist two independent intermediate edges between $I_a$ and $I_b$, which by Lemma 2 yields $|I_a| + |I_b| \geq 2\overline{p} + 8 = 10$, contradicting Claim 2. So, $N(w_1) \cap V(I_b^*) = \emptyset$. Further, if $\Upsilon(I_a, I_c) \neq \emptyset$ for some $c \in \{1, 2, ..., s\} \backslash \{a, b\}$ then by Lemma 2, $|I_a| + |I_c| \geq 2\overline{p} + 6 = 8$, implying that $|I_c| \geq 5$. But then $|I_b| + |I_c| \geq 10$, contradicting Claim 2. Hence $\Upsilon(I_a, I_i) = \emptyset$ for each $i \in \{1, 2, ..., s\} \backslash \{a, b\}$. Finally, if $w_1 \xi_{a+1} \in E(G)$ then
$$\xi_a x_1 x_2 \xi_b \overleftarrow{C} \xi_{a+1} w_1 w_2 w_4 \overrightarrow{C} \xi_a$$
is longer than $C$, a contradiction. Thus
$$N(w_1) \subseteq (\{\xi_1, \xi_2, ..., \xi_s\} \cup \{w_2\}) \backslash \{\xi_{a+1}\},$$
which yields $|N(w_1)| \leq s = \delta - 1$, a contradiction.

**Case 1.2.1.2.** $|I_a| = |I_b| = 4$.
Put $I_a = \xi_a w_1 w_2 w_3 \xi_{a+1}$ and $I_b = \xi_b w_4 w_5 w_6 \xi_{b+1}$.

**Case 1.2.1.2.1.** $y \in \{w_1, w_3\}$.

Assume w.l.o.g. that $y = w_3$. If $z \in \{w_5, w_6\}$ than $|C_2| > |C|$, a contradiction. Hence $z = w_4$.

**Case 1.2.1.2.1.1.** $\xi_{a+1} \neq \xi_b$.

If $w_1 \xi_{a+1} \in E(G)$ then

$$|\xi_a x_1 x_2 \xi_b \overleftarrow{C} \xi_{a+1} w_1 w_2 w_3 w_4 \overrightarrow{C} \xi_a| = |C| + 2,$$

a contradiction. Next, if $w_1 \xi_b \in E(G)$ then

$$|\xi_a x_1 x_2 \xi_{a+1} \overrightarrow{C} \xi_b w_1 w_2 w_3 w_4 \overrightarrow{C} \xi_a| = |C| + 2,$$

a contradiction. Further, if $w_1 w_4 \in E(G)$ then

$$|\xi_a x_1 x_2 \xi_b \overleftarrow{C} w_1 w_4 \overrightarrow{C} \xi_a| = |C| + 2,$$

again a contradiction. Moreover, if $N(w_1) \cap V(I_b^*) \neq \emptyset$ then there exist two independent intermediate edges between $I_a$ and $I_b$ which by Lemma 2 yields $|I_a| + |I_b| \geq 2\overline{p} + 8 \geq 10$, contradicting Claim 2. Furthermore, if $N(w_1) \cap V(I_i^*) = \emptyset$ for each $i \in \{1, 2, ..., s\} \backslash \{a, b\}$ then

$$N(w_1) \subseteq (\{\xi_1, \xi_2, ..., \xi_s\} \cup \{w_2, w_3\}) \backslash \{\xi_{a+1}, \xi_b\},$$

implying that $|N(w_1)| \leq s = \delta - 1$, a contradiction. Otherwise, $w_1 v \in E(G)$, where $v \in V(I_c^*)$ for some $c \in \{1, 2, ..., s\} \backslash \{a, b\}$. By a similar way, it can be shown that $w_2 u \in E(G)$, where $u \in V(I_d^*)$ for some $d \in \{1, 2, ..., s\} \backslash \{a, b\}$. By Lemma 2, $|I_a| + |I_c| \geq 2\overline{p} + 6 = 8$, that is $|I_c| \geq 4$. By Claim 2, $|I_c| = 4$. By a symmetric argument, $|I_d| = 4$. Put $I_c = \xi_c w_7 w_8 w_9 \xi_{c+1}$. As in Case 1.2.1.2.1, we can show that $v = w_9$, i.e. $w_1 w_9 \in E(G)$. If $d = c$ then $|\Upsilon(I_a, I_c)| = 2$ and by Lemma 2, $|I_a| + |I_c| \geq 2\overline{p} + 7 = 9$, a contradiction. Otherwise, there are at least four elementary segments of length at least 4, contradicting Claim 2.

**Case 1.2.1.2.1.2.** $\xi_{a+1} = \xi_b$.

Assume w.l.o.g. that $a = 1$ and $b = 2$. If $\Upsilon(I_1, I_2, ..., I_s) = \Upsilon(I_1, I_2) = \{w_3 w_4\}$ then by Claim 1, $\tau \leq 1$, a contradiction. Otherwise, there is an intermediate edge $uv$ such that $u \in V(I_1^*) \cup V(I_2^*)$ and $v \in V(I_c^*)$ for some $c \in \{1, 2, ..., s\} \backslash \{1, 2\}$. Assume w.l.o.g. that $u \in V(I_1^*)$. If $u = w_3$ then as above, $\xi_2 = \xi_c$, a contradiction. Let $u \neq w_3$. By Lemma 2, $|I_1| + |I_c| \geq 8$, i.e. $|I_c| \geq 4$. By Claim 2, $|I_c| = 4$. Put $I_c = \xi_c w_7 w_8 w_9 \xi_{c+1}$.

**Case 1.2.1.2.1.2.1.** $u = w_1$.

As in Case 1.2.1.2.1, $uv = w_1 w_9$. Then

$$|\xi_1 w_1 w_9 \overleftarrow{C} w_4 w_3 \xi_2 x_2 x_1 \xi_{c+1} \overrightarrow{C} \xi_1| \geq |C| + 1,$$

a contradiction.

**Case 1.2.1.2.1.2.2.** $u = w_2$.
If $v \in \{w_8, w_9\}$ then

$$|\xi_1 w_1 w_2 w_8 \overleftarrow{C} w_4 w_3 \xi_2 x_2 x_1 \xi_{c+1} \overrightarrow{C} \xi_1| \geq |C| + 1,$$

a contradiction. If $v = w_7$ then

$$|\xi_1 x_1 x_2 \xi_c \overleftarrow{C} w_2 w_7 \overrightarrow{C} \xi_1| = |C| + 1,$$

again a contradiction.

**Case 1.2.1.2.2.** $y = w_2$.
If $z = w_4$ then $|C_1| > |C|$, a contradiction. If $z = w_6$ then $|C_2| > |C|$, a contradiction. Hence $z = w_5$. Clearly, $\Upsilon(I_a, I_b) = \{w_2 w_5\}$. If $|I_i| = 3$ for each $i \in \{1, 2, ..., s\} \backslash \{a, b\}$ then by Lemma 2, $\Upsilon(I_1, I_2, ..., I_s) = \{w_2 w_5\}$ and by Claim 1, $\tau \leq 1$, contradicting the hypothesis. Otherwise, $|I_c| \geq 4$ for some $c \in \{1, 2, ..., s\} \backslash \{a, b\}$ and $|I_i| = 3$ for each $i \in \{1, 2, ..., s\} \backslash \{a, b, c\}$. By Claim 2, $|I_c| = 4$. Put $I_c = \xi_c w_7 w_8 w_9 \xi_{c+1}$. Clearly, $\Upsilon(I_1, I_2, ..., I_s) = \Upsilon(I_a, I_b, I_c)$. If $\Upsilon(I_a, I_c) = \Upsilon(I_b, I_c) = \emptyset$ then again $\tau \leq 1$, a contradiction. Let $uv \in E(G)$, where $u \in I_a^* \cup I_b^*$ and $v \in V(I_c^*)$. Assume w.l.o.g. that $u \in V(I_a^*)$. If $u \in \{w_1, w_3\}$ then we can argue as in Case 1.2.1.2.1. Let $u = w_2$, implying that $v = w_8$. Observing that $\{w_1, w_3, w_4, w_6, w_7, w_9\}$ is an independent set of vertices, we conclude that $G \backslash (\{\xi_1, \xi_2, ..., \xi_s\} \cup \{w_2, w_5, w_8\})$ has at least $s + 4$ connected components, that is $\tau < 1$, contradicting the hypothesis.

**Case 1.2.2.** $|I_a| + |I_b| = 9$.
Since $|I_i| \geq 3$ ($i = 1, 2, ..., s$), we can assume w.l.o.g. that either $|I_a| = 3$, $|I_b| = 6$ or $|I_a| = 4$, $|I_b| = 5$.

**Case 1.2.2.1.** $|I_a| = 3$ and $|I_b| = 6$.
By Claim 2, $|I_i| = 3$ for each $i \in \{1, 2, ..., s\} \backslash \{b\}$. Put

$$I_a = \xi_a w_1 w_2 \xi_{a+1}, \quad I_b = \xi_b w_3 w_4 w_5 w_6 w_7 \xi_{b+1}.$$

Since $|I_a| = 3$, we can assume w.l.o.g. that $y = w_2$. If $z = w_3$ then $|C_1| > |C|$, a contradiction. If $z \in \{w_6, w_7\}$ then $|C_2| > |C|$, a contradiction. So, $z \in \{w_4, w_5\}$.

**Case 1.2.2.1.1.** $z = w_4$.
If $w_1 w_4 \in E(G)$ then

$$|\xi_a x_1 x_2 \xi_b \overleftarrow{C} w_1 w_4 \overrightarrow{C} \xi_a| \geq |C| + 1,$$

a contradiction. Next, if $N(w_1) \cap V(I_b^*) \neq \emptyset$ then there are two independent intermediate edges between $I_a$ and $I_b$ and by Lemma 2, $|I_a| + |I_b| \geq 2\overline{p} + 8 = 10$, contradicting Claim 2. Further, if $w_1 \xi_{a+1} \in E(G)$ then

$$\xi_a x_1 x_2 \xi_b \overleftarrow{C} \xi_{a+1} w_1 w_2 w_4 \overrightarrow{C} \xi_a$$

is longer than $C$, a contradiction. Finally, by Claim 2, $N(w_1) \cap V(I_i^*) = \emptyset$ for each $i \in \{1, 2, ..., s\} \backslash \{a, b\}$. So,

$$N(w_1) \subseteq (\{\xi_1, \xi_2, ..., \xi_s\} \backslash \{\xi_{a+1}\}) \cup \{w_2\},$$

that is $|N(w_1)| \leq s = \delta - 1$, a contradiction.

**Case 1.2.2.1.2.** $z = w_5$.

If $w_2w_4 \in E(G)$ then we can argue as in Case 1.2.2.1.1. Let $w_2w_4 \notin E(G)$. It means that $w_5$ belongs to all intermediate edges. By Claim 1, $\tau \leq 1$, contradicting the hypothesis.

**Case 1.2.2.2.** $|I_a| = 4$ and $|I_b| = 5$.

By Claim 2, $|I_i| = 3$ and $\Upsilon(I_a, I_i) = \emptyset$ for each $i \in \{1, 2, ..., s\} \backslash \{a, b\}$. If $\Upsilon(I_b, I_c) \neq \emptyset$ for some $c \in \{1, 2, ..., s\} \backslash \{a, b\}$ then we can argue as in Case 1.2.1.1. Otherwise, $\Upsilon(I_1, I_2, ..., I_s) = \Upsilon(I_a, I_b)$. If There are two independent edges in $\Upsilon(I_a, I_b)$ then by Lemma 2, $|I_a| + |I_b| \geq 10$, contradicting Claim 2. Otherwise, by Claim 1, $\tau \leq 1$, a contradiction.

**Case 2.** $2 \leq \overline{p} \leq \delta - 3$.

It follows that $|N_C(x_i)| \geq \delta - \overline{p} \geq 3$ $(i = 1, 2)$. If $N_C(x_1) \neq N_C(x_2)$ then by Lemma 1, $|C| \geq 4\delta - 2\overline{p} \geq 3\delta - \overline{p} + 3$, contradicting (1). Hence $N_C(x_1) = N_C(x_2)$, implying that $|I_i| \geq \overline{p} + 2$ $(i = 1, 2, ..., s)$. Clearly, $s \geq |N_C(x_1)| - (|V(P)| - 1) \geq \delta - \overline{p} \geq 3$. If $s \geq \delta - \overline{p} + 1$ then

$$|C| \geq s(\overline{p} + 2) \geq (\delta - \overline{p} + 1)(\overline{p} + 2)$$

$$= (\delta - \overline{p} - 1)(\overline{p} - 1) + 3\delta - \overline{p} + 1 \geq 3\delta - \overline{p} + 3,$$

again contradicting (1). Hence $s = \delta - \overline{p}$. It means that $x_1x_2 \in E(G)$, that is $G[V(P)]$ is hamiltonian. By symmetric arguments, $N_C(y) = N_C(x_1)$ for each $y \in V(P)$. If $\Upsilon(I_1, I_2, ..., I_s) = \emptyset$ then by Claim 1, $\tau \leq 1$, contradicting the hypothesis. Otherwise $\Upsilon(I_a, I_b) \neq \emptyset$ for some elementary segments $I_a$ and $I_b$. By definition, there is an intermediate path $L$ between $I_a$ and $I_b$. If $|L| \geq 2$ then by lemma 2,

$$|I_a| + |I_b| \geq 2\overline{p} + 2|L| + 4 \geq 2\overline{p} + 8.$$

Hence

$$|C| = |I_a| + |I_b| + \sum_{i \in \{1, ..., s\} \backslash \{a, b\}} |I_i| \geq 2\overline{p} + 8 + (s - 2)(\overline{p} + 2)$$

$$= (\delta - \overline{p} - 2)(\overline{p} - 1) + 3\delta - \overline{p} + 2 \geq 3\delta - \overline{p} + 3,$$

contradicting (1). Thus, $|L| = 1$, i.e. $\Upsilon(I_1, I_2, ..., I_s) \subseteq E(G)$. By Lemma 2,

$$|I_a| + |I_b| \geq 2\overline{p} + 2|L| + 4 = 2\overline{p} + 6,$$

which yields

$$|C| = |I_a| + |I_b| + \sum_{i \in \{1,\dots,s\}\setminus\{a,b\}} |I_i| \geq 2\overline{p} + 6 + (s-2)(\overline{p}+2)$$

$$= (s-2)(\overline{p}-2) + (\delta - \overline{p} - 4) + (3\delta - \overline{p} + 2).$$

If either $\overline{p} \geq 3$ or $\delta - \overline{p} \geq 4$ then $|C| \geq 3\delta - \overline{p} + 2$, contradicting (1). Otherwise $\overline{p} = 2$ and $\delta - \overline{p} = 3$, implying that $s = 3$ and $\delta = 5$. Since $s = 3$, we can assume w.l.o.g. that $a = 1$ and $b = 2$, i.e. $|I_1| + |I_2| \geq 10$, $|I_3| \geq 4$ and $|C| \geq 14$. On the other hand, by (1), $|C| \leq 3\delta - \overline{p} + 1 = 14$, implying that

$$|I_1| + |I_2| = 10, \quad |I_3| = 4, \quad |C| = 14.$$

If $|I_1| = |I_2| = 5$ then by Lemma 2, $|\Upsilon(I_1, I_2)| \leq 1$ and $\Upsilon(I_1, I_3) = \Upsilon(I_2, I_3) = \emptyset$, implying that $|\Upsilon(I_1, I_2, I_3)| \leq 1$. By Claim 1, $\tau \leq 1$, contradicting the hypothesis. Now assume w.l.o.g. that $|I_1| = 4$ and $|I_2| = 6$. Let $L = yz$, where $y \in V(I_1^*)$ and $z \in V(I_2^*)$. Put

$$I_1 = \xi_1 w_1 w_2 w_3 \xi_2, \quad I_2 = \xi_2 w_4 w_5 w_6 w_7 w_8 \xi_3,$$

$$C_1 = \xi_1 x_1 \overrightarrow{P} x_2 \xi_2 \overleftarrow{C} yz \overrightarrow{C} \xi_1,$$

$$C_2 = \xi_1 \overrightarrow{C} yz \overleftarrow{C} \xi_2 x_1 \overrightarrow{P} x_2 \xi_3 \overrightarrow{C} \xi_1.$$

**Case 2.1**. $y = w_2$.

If $z \in \{w_4, w_5\}$ then $|C_1| > |C|$, a contradiction. Next, if $z \in \{w_7, w_8\}$ then $|C_2| > |C|$, again a contradiction. Hence, $z = w_6$. Further, if $w_1 w_3 \in E(G)$ then

$$\xi_1 x_1 \overrightarrow{P} x_2 \xi_2 w_3 w_1 w_2 w_6 \overrightarrow{C} \xi_1$$

is longer than $C$, a contradiction. Next, by Lemma 2, $N(w_1) \cap (V(I_2^*) \cup V(I_3^*)) = \emptyset$. Hence $N(w_1) \subseteq \{\xi_1, \xi_2, \xi_3, w_2\}$, that is $|N(w_1)| \leq 4$, contradicting the fact that $\delta = 5$.

**Case 2.2**. $y = w_1$.

If $z \in \{w_4, w_5, w_6\}$ then $|C_1| > |C|$, a contradiction. Next, if $z = w_8$ then $|C_2| > |C|$, a contradiction. Hence, $z = w_7$. Further, if $w_3 \xi_1 \in E(G)$ then

$$\xi_1 w_3 w_2 w_1 w_7 \overleftarrow{C} \xi_2 x_1 \overrightarrow{P} x_2 \xi_3 \overrightarrow{C} \xi_1$$

is longer than $C$, a contradiction. Observing also that by Lemma 2,

$$N(w_3) \cap (V(I_2^*) \cup V(I_3^*)) = \emptyset,$$

we get $N(w_3) \subseteq \{\xi_2, \xi_3, w_1, w_2\}$, that is $|N(w_3)| \leq 4$, a contradiction.

**Case 2.3**. $y = w_3$.

If $z = w_4$ then $|C_1| > |C|$, a contradiction. Next, if $z \in \{w_6, w_7, w_8\}$ then $|C_2| > |C|$, a contradiction. Hence, $z = w_5$. Further, if $w_1 \xi_2 \in E(G)$ then

$$\xi_1 x_1 \overrightarrow{P} x_2 \xi_2 w_1 w_2 w_3 w_5 \overrightarrow{C} \xi_1$$

is longer than $C$, a contradiction. Recalling also that $N(w_1) \cap (V(I_2^*) \cup V(I_3^*)) = \emptyset$, we get $N(w_1) \subseteq \{\xi_1, \xi_3, w_2, w_3\}$, contradicting the fact that $\delta = 5$.

**Case 3.** $2 \leq \overline{p} = \delta - 2$.

It follows that $|N_C(x_i)| \geq \delta - \overline{p} = 2$ ($i = 1, 2$). If $N_C(x_1) \neq N_C(x_2)$ then by Lemma 1, $|C| \geq 4\delta - 2\overline{p} = 3\delta - \overline{p} + 2$, contradicting (1). Hence, $N_C(x_1) = N_C(x_2)$. Clearly, $s = |N_C(x_1)| \geq 2$. Further, if $s \geq 3$ then

$$|C| \geq s(\overline{p} + 2) \geq 3\delta \geq 3\delta - \overline{p} + 2,$$

again contradicting (1). Hence, $s = 2$. It follows that $x_1 x_2 \in E(G)$, that is $G[V(P)]$ is hamiltonian. By symmetric arguments, $N_C(v) = N_C(x_1) = \{\xi_1, \xi_2\}$ for each $v \in V(P)$. If $\Upsilon(I_1, I_2) = \emptyset$ then by Claim 1, $\tau \leq 1$, contradicting the hypothesis. Otherwise, there is an intermediate path $L = yz$ such that $y \in V(I_1^*)$ and $z \in V(I_2^*)$. If $|L| \geq 2$ then by Lemma 2,

$$|C| = |I_1| + |I_2| \geq 2\overline{p} + 2|L| + 4 \geq 2\overline{p} + 8 = 3\delta - \overline{p} + 2,$$

contradicting (1). hence $|L| = 1$, implying that $\Upsilon(I_1, I_2) \subseteq E(G)$. If there are two independent intermediate edges between $I_1, I_2$, then by Lemma 2, $|C| = |I_1| + |I_2| \geq 2\overline{p} + 8 = 3\delta - \overline{p} + 2$, contradicting (1). Otherwise, by Claim 1, $\tau \leq 1$, contradicting the hypothesis.

**Case 4.** $2 \leq \overline{p} = \delta - 1$.

It follows that $|N_C(x_i)| \geq \delta - \overline{p} = 1$ ($i = 1, 2$). By (1), $3\delta - 3 \leq |C| \leq 2\delta + 2$, implying that $\delta \leq 5$

**Case 4.1.** $|N_C(x_i)| \geq 2$ ($i = 1, 2$).

If $N_C(x_1) \neq N_C(x_2)$ then by Lemma 1, $|C| \geq 2\overline{p} + 8 = 3\delta - \overline{p} + 5$, contradicting (1). Hence, $N_C(x_1) = N_C(x_2)$. Clearly $s \geq 2$. Further, if $s \geq 3$ then

$$|C| \geq s(\overline{p} + 2) \geq 3(\delta + 1) > 3\delta - \overline{p} + 2,$$

contradicting (1). Let $s = 2$. It follows that $|C| \geq s(\overline{p} + 2) = 2\delta + 2$. By (1), $|C| \leq 2\delta + 2$, implying that

$$|C| = 2\delta + 2, \ |I_1| = |I_2| = \delta + 1, \ V(G) = V(C \cup P).$$

Since $\kappa \geq 3$, there is an edge $zw$ such that $z \in V(P)$ and $w \in V(C) \backslash \{\xi_1, \xi_2\}$. Assume w.l.o.g. that $w \in V(I_1^*)$. Then it is easy to see that $|I_1| \geq \delta + 3$, a contradiction.

**Case 4.1.** Either $|N_C(x_1)| = 1$ or $|N_C(x_2)| = 1$.

Assume w.l.o.g. that $|N_C(x_1)| = 1$. Put $N_C(x_1) = \{y_1\}$. If $N_C(x_1) \neq N_C(x_2)$ then $x_2y_2 \in E(G)$ for some $y_2 \in V(C)\backslash\{y_1\}$ and we can argue as in Case 4.1. Let $N_C(x_1) = N_C(x_2) = \{y_1\}$. Since $\kappa \geq 1$, there is an edge $zw$ such that $z \in V(P)$ and $w \in V(C)\backslash\{y_1\}$. Clearly, $z \notin \{x_1, x_2\}$ and $x_2z^- \in E(G)$, where $z^-$ is the previous vertex of $z$ along $\overrightarrow{P}$. Then replacing $P$ with $x_1\overrightarrow{P}z^-x_2\overleftarrow{P}z$, we can argue as in Case 4.1.

**Case 5**. $\delta \leq \overline{p} \leq \delta + 1$.

By (1), $|C| \leq 2\delta + 1$. By Theorem D, $C$ is a Hamilton cycle, contradicting the fact that $\overline{p} \geq \delta \geq 3$. ∎

Institute for Informatics and Automation Problems
National Academy of Sciences
P. Sevak 1, Yerevan 0014, Armenia
E-mail: zhora@ipia.sci.am